\documentclass[a4paper]{article}
\usepackage{amsmath, amssymb, eucal}

\newtheorem{thm}{Theorem}[section]
\newtheorem{lem}[thm]{Lemma}
\newtheorem{eg}[thm]{Example}
\newtheorem{prop}[thm]{Proposition}
\newtheorem{cor}[thm]{Corollary}
\newtheorem{defn}[thm]{Definition}
\newtheorem{lem-defn}[thm]{Lemma and Definition}
\newtheorem{rem}[thm]{Remark}
\newtheorem{ntn}[thm]{Notation}
\newtheorem{rem-eg}[thm]{Remark and Example}

\newcommand{\smnoind}{{\smallskip\noindent}}

\newcommand{\id}{{\rm id}}

\newcommand{\ti}{\tilde}

\newcommand{\abs}[1]{\left\vert#1\right\vert}

\newcommand{\real}{{\rm Re}}

\newcommand{\h}{\mathbb{H}}
\newcommand{\R}{\mathbb{R}}
\newcommand{\C}{\mathbb{C}}
\newcommand{\re}{\mathbf{Re}}
\newcommand{\ep}{\epsilon}
\newcommand{\hthr}[1]{\h \otimes_{#1}\ \!_\triangleright \h _\triangleleft}

\newcommand{\hthlr}[1]{\ \!_\triangleright \h \otimes_{#1} \h _\triangleleft}
\newcommand{\mcl}[3]{\mathcal{L}_{#1}(#2;#3)}
\newcommand{\ml}[2]{\mathcal{L}_{#1}(#2)}
\newcommand{\cl}{\mathcal{L}}
\newcommand{\lala}{\langle\!\langle}
\newcommand{\rara}{\rangle\!\rangle}
\newcommand{\la}{\langle}
\newcommand{\ra}{\rangle}

\newenvironment{prf}{{\noindent \textbf{Proof:}\ }}{\hfill $\Box$\\ \smallskip}

\begin{document}

\title{On quaternionic functional analysis}
\author{Chi-Keung Ng\thanks{This work is supported by the National Natural Science Foundation of China (10371058).}}
\date{}
\maketitle
\begin{abstract}
In this article, we will show that the category of quaternion vector spaces, the category of (both one-sided and two sided) quaternion Hilbert spaces and the category of quaternion $B^*$-algebras are equivalent to the category of real vector spaces, the category of real Hilbert spaces and the category of real $C^*$-algebras respectively. 
We will also give a Riesz representation theorem for quaternion Hilbert spaces and will extend the main results in \cite{Kul} (namely, we will give the full versions of the Gelfand-Naimark theorem and the Gelfand theorem for quaternion $B^*$-algebras). 
On our way to these results, we compare, clarify and unify the term ``quaternion Hilbert spaces'' in the literatures.

\end{abstract}

\bigskip

\section*{0\hspace{0.225in}Introduction}

\bigskip

Recently, there are some interests in ``quaternion Hilbert spaces'' (or Wachs spaces), ``quaternion normed spaces'' as well as ``quaternion algebras'' (see \cite{AJ}, \cite{AK}, \cite{AK2}, \cite{BCT}, \cite{CT}, \cite{BH}, \cite{HR-ten}, \cite{HR-proj}, \cite{HR-uniq}, \cite{HS}, \cite{Kul}, \cite{Lud}, \cite{NV}, \cite{Tor0}, \cite{Tor}, \cite{Tor2}, \cite{Tor-dual} \& \cite{Vis}).
The aim of this article is to give a systematic study of these objects. 
Note that although the ``quaternion Hilbert spaces'' considered in \cite{HS}, \cite{Tor} and \cite{Vis} are one-sided, they are automatically two-sided (see e.g. Lemma \ref{cp-hil}). 
The starting point of this research is an observation we found in the argument of \cite[Theorem 4.1]{HS}.
We believe that one can use this observation to obtain a lot of results in ``quaternionic functional analysis'' (which is related to ``quaternionic quantum mechanics''; see e.g. \cite{CT}, \cite{BH}, \cite{HR-ten}, \cite{HR-proj}, \cite{HR-uniq} \& \cite{HS}). 
As demonstrations, we will use this observation to generalise and to simplify the main results in \cite{HS} and \cite{Kul} as well as to give the results stated in the abstract. 

\medskip

We will start with (two-sided) quaternion vector spaces and show in Section 1 that the category of quaternion vector spaces is isomorphic to the category of real vector spaces. 
Based on this, we show in Section 2 that every (two-sided) quaternion normed space is naturally a ``quaternionization'' of a real normed space. 
However, since there can be many different quaternionizations of a given real normed space, the category of quaternion normed spaces is in general ``larger''. 
We then show that if $X$ is a quaternion normed space, then $\mcl{\h}{X}{\h} \cong \mcl{\R}{X_\real}{\R}$ (as real Banach spaces) where $X_\real$ is the space of ``centralizers'' for the quaternion scalar multiplication. 
Using this, we obtained very easily \cite[Corollary 4.1]{HS} as well as the Hahn Banach theorem as stated in \cite[Theorem 4.1]{HS} (in fact, we have formal extensions of these two results).

\medskip

In Section 3, we will give three notions of ``quaternion Hilbert spaces''. 
We will show that they are all the same and are (categorically) equivalent to real Hilbert spaces.  
Moreover, we will give a Riesz representation type theorem for them. 

\medskip

Finally, we will show in Section 4 that the category of quaternion $B^*$-algebras is equivalent to the category of real $C^*$-algebras. 
On our way to this equivalence, we will give, using very simple arguments, two versions of the Gelfand-Naimark theorem as well as the Gelfand theorem for quaternion $B^*$-algebras (which improve the corresponding results in \cite{Kul}). 

\bigskip

\bigskip

\section{Quaternion vector spaces}

\bigskip

\begin{defn}
Let $D$ be a division ring and $F$ be the center of $D$.
A unital $D$-bimodule $X$ is called a \emph{$D$-vector space} if $X$ is a $F$-vector space under the restriction of the scalar multiplication to $F$.
\end{defn}

\medskip

In this article, we mainly concern with the case when $D = \h$, the set of all real quaternions.
Let us recall its definition here. 
Suppose that $\h$ is a four dimensional $\R$-vector space with basis $B := \{1,i,j,k\}$ and consider a $\R$-linear map from $\h$ to $M_4(\R)$ given by
\begin{equation}
\label{h<M4}
a1+bi+cj+dk \longmapsto \left(\begin{array}{cccc}
a  & b  & c  & d\\
-b & a  & -d & c\\
-c & d  & a  & -b\\
-d & -c & b  & a
\end{array}\right).
\end{equation}
This map induces an involutive algebra structure on $\h$ such that $1$ is the identity and
$$i^2\ =\ j^2\ =\ k^2\ =\ -1, \quad ij\ =\ k\ =\ -ji, \quad jk\ =\ i\ =\ -kj, \quad ki\ = \ j\ =\ -ik.$$ 
Moreover, the conjugation is given by 
$$(a1+bi+cj+dk)^*\ =\ a1 -bi-cj-dk.$$
In the following we will write $a+bi+cj+dk$ instead of $a1+bi+cj+dk$.

\medskip

\begin{rem}
\label{unital}
\rm
The notion of $\h$-vector space appears in many literatures in disguised forms. 
The following are two examples.

\smnoind
(a) In \cite{HS}, Horwitz and Soffer considered a real $B^*$-algebra $A$ containing a real $*$-subalgebra that is $*$-isomorphic to $\h$. 
However, it seems that they implicitly assume $A$ being a real unital $B^*$-algebra containing a unital real $*$-subalgebra $B$ which is $*$-isomorphic to $\h$.
This assumption is definitely needed in \cite[Corollary 4.1]{HS}:
for any quaternion linear functional $\rho$ and any $a\in A$, one has 
$\rho(A_1a) = 1\rho(a) = \rho(a)$ (where $A_1\in B$ is the corresponding element of $1\in \h$) and by \cite[Corollary 4.1]{HS}, one has $A_1 a = a$; similarly, one has $a A_1 = a$. 
With the above unital assumption, those real $B^*$-algebras considered in \cite{HS} are actually $\h$-vector spaces. 

\smnoind
(b) Although in \cite[Example 2]{Kul}, $\h$ is not a unital subalgebra of $A$ but in \cite[Theorem 4]{Kul}, one needs to assume that $\h$ is a unital subalgebra of $A$ (see Remark \ref{rem-hil}(a) below). 
Therefore, the algebra concerned in the main result of \cite{Kul} is again a $\h$-vector space. 
\end{rem}

\medskip

\begin{ntn}
\rm 
Throughout this article, we will identify $\R$ with the center of $\h$ and  $\otimes$ means the algebraic tensor product over $\R$.
\end{ntn}

\medskip

\begin{eg}
(a) It is clear that $\h$ is itself a $\h$-vector space.

\smnoind
(b) $\h \otimes \h$ together with the scalar multiplication:
$$\alpha (a\otimes b) \beta = a \otimes \alpha b \beta$$
is a $\h$-vector space and is denoted by $\hthr{\!}$.
We can define, in a similar fashion, $\ \!_\triangleright \h _\triangleleft \otimes\h$. 
These two $\h$-vector spaces are regarded as ``of the same type'' as they differ only by an ``external swapping''. 

\smnoind
(c) One can define another type of scalar product on $\h \otimes \h$ by:
$$\alpha (a\otimes b) \beta = \alpha a \otimes b \beta$$
and we denote this $\h$-vector space by $\hthlr{}$. 
Similarly, we can define $\h _\triangleleft \otimes\ \!_\triangleright \h$ (which is ``of the same type'' as $\hthlr{}$). 
\end{eg}

\medskip

It is not easy to see directly whether $\hthlr{}$ is isomorphic to $\hthr{\!}$ as $\h$-vector spaces (although one can show that this is the case using Theorem \ref{qvs=rvs} below). 
Note that $1\otimes 1$ ``generates'' $\hthlr{}$ but it is not an obvious task to find a generator for $\hthr{\!}$.

\medskip

The starting point of this article is the following somewhat surprising ``polarization result'' which is hidden in the argument of \cite[Theorem 4.1]{HS}.
Since this result follows from direct computation, its proof will be omitted. 

\medskip

\begin{lem}
\label{polar-decom}
Let $X$ be a $\h$-vector space and $X_\real := \{ x\in X: \alpha x = x \alpha {\rm\ for\ any\ } \alpha\in \h\}$ (called the \emph{real part} of $X$).
Define a $\R$-linear map $\re : X \rightarrow X$ by
$$\re(x)\ :=\ \frac{1}{4} \sum_{e\in B} e^* xe \qquad (x\in X).$$
Then $\re(X) = X_\real$, $\re\circ \re = \re$ and for any $x\in X$, we have the \emph{polarization identity}:
\begin{equation}
\label{pol-id}
x\ =\ \sum_{e\in B} \re( e^* x)e.
\end{equation}
\end{lem}

\medskip

\begin{rem}
\label{rem-pol}
\rm
(a) Note that the decomposition in (\ref{pol-id}) is unique in the sense that if $x\in X$ and $x_e\in X_\real$ ($e\in B$) such that $x = \sum_{e\in B} x_e e$, then $x_e = \re( e^* x)$ for any $e\in B$. 

\smnoind
(b) $X_\real$ is the unique real vector space (up to isomorphism) whose \emph{quaternionization}, $X_\real\otimes \h$, is isomorphic to $X$ as $\h$-vector space. 
\end{rem}

\medskip

In the following, we denote by $L_\h(X;Y)$ (respectively, $L_\R(V;W)$) the space of all $\h$-bimodule (respectively, $\R$-linear) maps from a $\h$-vector space $X$ (respectively, $\R$-vector space $V$) to another $\h$-vector space $Y$ (respectively, $\R$-vector space $W$).

\medskip

\begin{lem}
\label{qmor=rmor}
Let $X$ and $Y$ be two $\h$-vector spaces.
If $T\in L_\h(X;Y)$, then $T\! \mid \!_{X_\real}\in L_\R(X_\real;Y_\real)$.
This induces a $\R$-linear isomorphism $\Psi$ from $L_\h(X;Y)$ onto $L_\R(X_\real;Y_\real)$.
\end{lem}
\begin{prf}
It is clear that $T(X_\real) \subseteq Y_\real$ and the map $\Psi: T\mapsto T\! \mid \!_{X_\real}$ is a $\R$-linear map from $L_\h(X;Y)$ to $L_\R(X_\real;Y_\real)$.
If $T\! \mid \!_{X_\real} = 0$, then the polarization identity (\ref{pol-id}) tells us that $T = 0$.
On the other hand, for any $S \in L_\R(X_\real;Y_\real)$, it is easy to see that the map 
$$\ti S \ =\ S\otimes \id_\h: X_\real \otimes \h \longrightarrow Y_\real\otimes \h $$ 
is a $\h$-bimodule map.
Under the identification of $X_\real \otimes \h $ with $X$ (Remark \ref{rem-pol}(b)), we have $\Psi(\ti S) = S$. 
\end{prf}

\medskip

The above discussions imply the following result. 
This result is a bit surprising when compare to the situation between real and complex vector spaces. 
Note that although a complex vector space is always a complexification of a real vector space, the choice of such a real vector space is not unique (but depends on the choice of a conjugation on that complex vector space) and thus the category of real vector spaces is not equivalent to that of complex vector spaces. 

\medskip

\begin{thm}
\label{qvs=rvs}
The category of $\h$-vector spaces is equivalent to the category of $\R$-vector spaces.
\end{thm}

\medskip

In a sense, this theorem ``trivializes'' the seemingly mysterious objects: ``quaternion vector spaces''. 
Furthermore, using this result, one can show easily that $\hthlr{} \cong \hthr{\!}$ as $\h$-vector spaces (the details will be given in \cite{Ng-dual-quat}). 

\medskip

We end this section with the following complexification result. 
Suppose that $\alpha_e \in \mathbb{R}$ ($ e \in \{i,j,k\}$) such that $\sum_{e\in \{i,j,k\}} \alpha_e^2 =1$. 
If $\alpha : = \sum_{e\in \{i,j,k\}} \alpha_e e \in \h$, then $\alpha^2 = -1$, $\abs{\alpha} =1$ and the $\mathbb{R}$-linear map $\Phi_\alpha: \mathbb{C} \rightarrow \mathbb{R} + \alpha\mathbb{R} \subseteq \h$ defined by $\Phi_\alpha(1) = 1$ and $\Phi_\alpha(i) = \alpha$
is an isometric algebra isomorphism. 
It is not hard to check that for any element $x = \sum_{e\in B} x_e e$ in a $\h$-vector space $X$, the equality $\alpha x = x \alpha$ is equivalent to the equalities $\alpha_e x_f = \alpha_f x_e$ (for any $e,f\in \{i,j,k\}$) and this is the case if and only if $x\in X_\real + \alpha X_\real$. 

\medskip

\begin{prop}
\label{complex}
Let $X$ be a $\h$-vector space. 
If $\alpha$ and $\Phi_\alpha$ are as in the above, then $X_{1,\alpha} := \{x\in X: \alpha x = x \alpha\}$ is a $\mathbb{C}$-vector space under the scalar multiplication induced by $\Phi_\alpha$ and we have $X_{1,\alpha} = X_\real + \alpha X_\real$. 
\end{prop}

\medskip

The above applies, in particular, to the case when $\alpha =i, j$ or $k$. 

\bigskip

\bigskip

\section{Quaternion normed spaces}

\bigskip

Although quaternion vector spaces are nothing more than real vector spaces, the notion of quaternion normed spaces still require some studies.
The reason is that there are many different ways to ``quaternionize'' a real normed space (similar to situation of complexifications of real normed spaces). 
In other words, the correspondence from quaternion normed spaces to real normed spaces is ``non-injective''
(in fact, one can show that given a real normed space $Z$, both the injective tensor norm and the projective tensor norm on $Z\otimes \h$ are ``norm quaternionizations'' of $Z$; see e.g \cite{Ng-dual-quat} for the details).

\medskip

Let us begin this section by recalling the norm on $\h$. 
If one identifies $M_4(\R)$ with $\ml{\R}{l^2_{(4)}(\R)}$ (where $l^2_{(4)}(\R)$ is the 4-dimensional real Hilbert space), then the norm induced on $\h$ through the map in (\ref{h<M4}) is:
$$\big\|a+bi+cj+dk\big\| \ := \ \sqrt{a^2 + b^2 + c^2 + d^2}$$
(therefore, we regarded $\h$ as a $\R$-$B^*$-subalgebra of $\ml{\R}{l^2_{(4)}(\R)}$). 
For any $\alpha\in \h\setminus\{0\}$, we have $\alpha^{-1} = \frac{ \alpha^*}{\|\alpha\|^2}$ and so, $\|\alpha^{-1}\| = \|\alpha\|^{-1}$.

\medskip

Note that $\h$ is a $\R$-Hilbert space and so $\h^*\cong \h$ as $\R$-normed spaces (where $\h^*$ is the dual space of $\h$). 
More precisely, for any $e\in B$, we can consider $\hat e\in \h^*$ as defined by $\hat e(f) = \delta_{e,f}$ ($f\in B$) and this extends by linearity to an isometry from $\h$ to $\h^*$. 
It is not hard to check that
$$\beta^*\cdot \hat\alpha \ =\ \widehat{\alpha\beta}\ = \ \hat \beta \cdot \alpha^* \qquad (\alpha, \beta\in \h).$$ 

\medskip

\begin{defn}
\label{def-q-norm}
A $\h$-vector space $X$ is called a \emph{$\h$-normed space} if there exists a norm on $X$ that turns it into a normed $\h$-bimodule (under the given $\h$-scalar multiplication).
A $\h$-normed space is called a \emph{$\h$-Banach space} if it is complete under the given norm. 
\end{defn}

\medskip

In the following, we denote by $\mcl{\h}{X}{Y}$ (respectively, $\mcl{\R}{V}{W}$) the space of all bounded $\h$-bimodule maps (respectively, bounded $\R$-linear maps) from a $\h$-normed space $X$ (respectively, $\R$-normed space $V$) to another $\h$-normed space $Y$ (respectively, $\R$-normed space $W$).
Moreover, we denote $\mcl{\R}{V}{\R}$ by $V^*$. 

\begin{rem}
\label{rem-dual}
\rm
Let $X$ be a $\h$-normed space. 

\smnoind
(a) One automatically has $\|\alpha x \beta\| = \|\alpha\| \ \|x\| \ \|\beta\|$ for any $x\in X$ and $\alpha, \beta\in \h$.
Indeed, if $\alpha\in \h \setminus \{0\}$, then
$$\|\alpha x\| \ \leq \ \|\alpha\| \ \|x\| \ \leq \ \|\alpha\| \ \|\alpha^{-1}\| \ \|\alpha x\| \ = \ \|\alpha\| \ \|\alpha\|^{-1} \ \|\alpha x\| \ = \ \|\alpha x\|.$$
Similarly, $\|x \beta \| = \|x\| \ \|\beta\|$.

\smnoind
(b) It is well known that $X^*$ is also a normed $\h$-bimodule (and hence a $\h$-normed space) under the canonical multiplication: $(\alpha \cdot f)(x) = f(x\alpha)$ and $(f\cdot \alpha)(x) = f(\alpha x)$ ($\alpha\in \h$; $x\in X$; $f\in X^*$). 
In this case, one has $(X^*)_\real \cong (X_\real)^*$ as $\R$-normed spaces. 
In fact, it is clear that $f\circ \re \in (X^*)_\real$ for any $f\in (X_\real)^*$ (because $\re(\alpha x) = \re(x\alpha)$). 
Conversely, for any $g\in (X^*)_\real$ and any $z\in X_\real$, we have 
$$g(iz)\ =\ (g\cdot j)(kz) \ = \ (j\cdot g)(kz) \ =\ -g(iz)$$
and similarly, $g(jz) = 0 = g(kz)$. 
Therefore, $g = (g\!\mid\!_{X_\real})\circ \re$ and it is easy to check that $f\mapsto f\circ \re$ is an isometry from 
$(X_\real)^*$ onto $(X^*)_\real$. 
Thus, from now on, we will write $X_\real^*$. 
\end{rem}

\medskip

The idea of the following lemma also comes from the proof of \cite[Theorem 4.1]{HS}.

\medskip

\begin{lem}
\label{lh=lr}
Let $X$ be a $\h$-normed space and $\re$ be the projection as defined in Lemma \ref{polar-decom}.

\smnoind
(a) $\|\re\| \leq 1$;

\smnoind
(b) $\mcl{\h}{X}{\h} \cong \mcl{\R}{X_\real}{\R} = X_\real^*$ as $\R$-Banach spaces (under the correspondence $\Psi$ as defined in Lemma \ref{qmor=rmor}).
\end{lem}
\begin{prf}
(a) This part is clear.

\smnoind
(b) It is clear that $\Psi: \mcl{\h}{X}{\h} \rightarrow X_\real^*$ is a $\R$-linear contraction.
For any $f\in X_\real^*$, define
$\ti f(x) = \sum_{e\in B} f(\re( e^* x))e$ (for $x\in X$).
It is easy to see that $\ti f\in \mcl{\h}{X}{\h}$ and $\Psi(\ti f) = f$.
If $x\in X$ with $\ti f(x) \neq 0$, we let $\alpha = \frac{\ti f(x)^*}{\|\ti f(x)\|}$ and
$\|\ti f(x) \| = \alpha \ti f(x) = \ti f (\alpha x) = f(\re(\alpha x))$ (as $\ti f (\alpha x) = \|\ti f(x) \| \in \R$ and we have equation (\ref{pol-id})). 
Therefore, $\|\ti f\| \leq \|f\|$ (by part (a)).
\end{prf}

\medskip

Using the above lemma, one can easily obtain the following result which was stated without proof in \cite[Corollary 4.1]{HS}. 
In the introduction of \cite{Kul}, Kulkarni raised an objection to this corollary with an example. 
However, in \cite[Example 2]{Kul}, $\h$ is not a unital subalgebra of the algebra $A$ and only \emph{positive} linear quaternion functionals are considered in that example. 
Therefore, Kulkarni's example does not contradict the following result (which is a clarified version of \cite[Corollary 4.1]{HS}). 

\medskip

\begin{cor}(\cite[Corollary 4.1]{HS})
\label{sep-pts}
$\mcl{\h}{X}{\h}$ separates points of $X$. 
\end{cor}
\begin{prf}
Suppose that $x\in X$ such that $T(x) = 0$ for any $T\in \mcl{\h}{X}{\h}$. 
Let $x = \sum_{e\in B} x_ee$ where $x_e\in X_\real$. 
Then $0 = T(x) = \sum_{e\in B} \Psi(T)(x_e)e$ which implies that $\Psi(T)(x_e) = 0$ for any $e\in B$. 
Since $T$ is arbitrary, Lemma \ref{lh=lr}(b) shows that $x_e = 0$ for any $e\in B$ and thus, $x=0$. 
\end{prf}

\medskip

There are different forms of Hahn Banach type theorem for ``$\h$-normed spaces'' in the literatures (depending on different settings and definitions; see e.g. \cite{Suh}, \cite{Tor0} and \cite{Tor2}). 
The following form of Hahn Banach theorem was proved in a slightly restrictive situation in \cite[Theorem 4.1]{HS}.
This theorem is actually an easy corollary of Lemma \ref{lh=lr}(b) as well as the Hahn Banach theorem for $\R$-normed spaces (note that since our Lemma \ref{lh=lr}(b) comes from the argument of \cite[Theorem 4.1]{HS}, we do not claim to give a new proof for this result).

\medskip

\begin{thm}(Horwitz-Razon)
\label{hb-h}
Suppose that $X$ is a $\h$-normed space and $Y$ is a $\h$-vector subspace of $X$.
For any $g\in \mcl{\h}{Y}{\h}$, there exists $\bar g\in \mcl{\h}{X}{\h}$ such that $\bar g\!\mid\!_Y = g$ and $\|\bar g\| = \| g\|$.
\end{thm}

\bigskip

\section{Quaternion Hilbert spaces}

\bigskip

In this section, we will discuss the notion of ``quaternion Hilbert spaces''. 
In fact, there are different definitions for this terminology in the literatures and we will show that they are essentially the same. 

\medskip

\begin{defn}
\label{def-hil}
(a) Let $Y$ be a right $\h$-module. 
A map $\langle \cdot, \cdot \rangle: Y \times Y \rightarrow \h$ is called a \emph{$\h$-valued inner product} on $Y$ if it satisfies the following conditions (for any $x,y,z\in Y$ and $\alpha, \beta\in \h$):

\begin{enumerate}

\item[i.] $\langle x, y+z\beta \rangle = \langle x, y \rangle + \langle x, z \rangle \beta$;

\item[ii.] $\langle y, x \rangle = \langle x, y \rangle^*$;

\item[iii.] $\langle x, x \rangle \geq 0$ and $\langle x, x \rangle = 0$ will imply that $x=0$. 

\end{enumerate}
In this case, $Y$ is called an \emph{inner product right $\h$-module}. 
Moreover, $Y$ is called a \emph{Hilbert right $\h$-module} (or a \emph{right Wachs space}; see e.g. \cite{Tor-dual}) if $Y$ is complete under the norm given by $\|x\| = \sqrt{\|\langle x, x \rangle\|}$. 

\smnoind
(b) An inner product right $\h$-module $Y$ is called an \emph{inner product $\h$-bimodule} if there exists a left $\h$-module structure on $Y$ such that $Y$ becomes a $\h$-vector space and
\begin{equation}
\label{compat}
\langle \alpha x, y \rangle \ = \ \langle x, \alpha^* y \rangle \qquad (x,y\in Y; \alpha\in \h).
\end{equation}
If $Y$ is complete with respect to the norm defined by $\|y\| := \sqrt{\|\langle y, y \rangle\|}$, then $Y$ is called a \emph{Hilbert $\h$-bimodule}. 

\smnoind
(c) Let 
$$(\h\otimes \h)_p\ :=\ \left\{\sum_{i=1}^n \alpha^*_i\otimes \alpha_i: n\in \mathbb{N}; \alpha_i\in \h\right\}$$ 
and  
$$(\alpha\otimes \beta)^\# \ :=\ \beta^*\otimes \alpha^* \qquad (\alpha, \beta\in \h).$$
A $\h$-vector space $Y$ is called a \emph{two-sided $\h$-inner product space} if there exists a map 
$$\lala \cdot, \cdot \rara: Y\times Y \rightarrow \h \otimes \h,$$ 
called a \emph{two sided $\h$-inner product}, satisfying the following conditions (for any $x,y,z\in Y$ and $\alpha, \beta\in \h$): 

\begin{enumerate}
\item[i.] $\lala x, \alpha y + z\beta \rara \ =\ (1\otimes \alpha)\lala x,y \rara + \lala x, z \rara (1\otimes \beta)$; 

\item[ii.] $\lala y,x \rara \ = \ \lala x, y \rara ^\#$; 

\item[iii.] $\lala x, x \rara \in (\h\otimes \h)_p$ and $\lala x, x \rara = 0$ if and only if $x = 0$.
\end{enumerate}
In this case, we can define a norm $\|x\| = \sqrt{m(\lala x, x \rara)}$ on $Y$ (where $m: \h \otimes \h \rightarrow \h$ is the multiplication) and we called $Y$ a \emph{two-sided $\h$-Hilbert space} if it is complete with respect to $\|\cdot\|$. 
\end{defn}

\medskip

\begin{rem}
\label{rem-hil}
\rm
(a) Suppose that $(Y, \langle \cdot, \cdot \rangle)$ is an inner product right $\h$-module. 
For any $x,y\in Y$, we have $\langle x, y\cdot 1\rangle = \langle x, y \rangle$ and so $\langle x, y - y\cdot 1 \rangle =0$. 
By putting $x = y - y\cdot 1 $, we see that $y\cdot 1 = y$ and so $Y$ is automatically a unital right $\h$-module. 
Similarly, a quaternionic Hilbert space as defined in \cite[Definition 3]{Kul} is a unital left $\h$-module. 

\smnoind
(b) It is easy to check that if $(Y, \lala \cdot, \cdot \rara)$ is a two-sided $\h$-inner product space, then $\lala \alpha y \beta, x \rara \ =\ (\beta^* \otimes 1)\ \lala y, x \rara  (\alpha^* \otimes 1)$ (for any $x,y\in Y$ and $\alpha, \beta\in \h$). 

\smnoind
(c) Note that in \cite[Theorem 4]{Kul}, one needs to consider \emph{unital} $\R$-$B^*$-algebras containing a \emph{unital} real $*$-subalgebra that is $*$-isomorphic to $\h$.
In fact, this assumption is needed in order to ensure that the left $\h$-module structure on the $\R$-Hilbert space $X$ in \cite[Theorem 4]{Kul} is unital (see part (a)) because $\pi(1)$ need to be the identity of $BL(X)$. 

\smnoind
(d) Let $(X, \langle \cdot, \cdot \rangle)$ be a Hilbert right $\h$-module.
Then the set $\ml{\h,r}{X}$ of all adjoinable (right) $\h$-module maps on $X$ is a $\R$-$B^*$-algebra (the argument for this fact is similar to its complex counterpart; see e.g. \cite{Lan}) and $X$ is unital because of part (a).
The existence of a left scalar multiplication on $X$ that turns $X$ into a $\h$-vector space
is the same as the existence of a unital homomorphism $\pi$ from $\h$ to $\ml{\h,r}{X}$. 
Furthermore, $\pi$ is a $*$-homomorphism if and only if the corresponding left scalar multiplication satisfies equality (\ref{compat}) (i.e. $X$ is a Hilbert $\h$-bimodule under this left scalar multiplication). 
\end{rem}

\medskip

Remark \ref{rem-hil}(d) tells us that Hilbert $\h$-bimodules are precisely \emph{essential $C^*$-correspondence over (the real $C^*$-algebra) $\h$}. 

\medskip

\begin{eg}
\label{eg-hil}
Let $K$ be a $\R$-inner product space. 
Then one can define a $\h$-valued inner product on $K\otimes \h$ in the canonical way (i.e. 
$\langle \xi\otimes \alpha, \eta\otimes \beta \rangle = \langle \xi, \eta \rangle\ \alpha^*\beta$) and this turns $K\otimes \h$ into a inner product $\h$-bimodule. 
Moreover, $K\otimes \h$ is a Hilbert $\h$-bimodule if $K$ is complete. 
On the other hand, one can also make $K\otimes \h$ into a two-sided $\h$-inner product space (respectively, $\h$-Hilbert space) in a canonical way (i.e. $\lala \xi\otimes \alpha, \eta\otimes \beta \rara = \langle \xi, \eta \rangle\ \alpha^*\otimes \beta$). 
\end{eg}

\medskip

The following is a known result. 
It contains a claim from \cite{HS} as well as several claims from \cite{Tor} and \cite{Vis}. 
Since no proof was found in those papers for these claims and they are crucial to our study, we will provide their simple arguments here. 

\medskip

\begin{prop}
\label{cp-hil}
Let $(X, \langle \cdot, \cdot \rangle)$ be a Hilbert right $\h$-module.

\smnoind
(a) There exists an orthonormal basis for $X$. 
Moreover, any two such bases are of the same cardinality.

\smnoind
(b) There exists a unital $\R$-$*$-homomorphism $\Lambda: \h \rightarrow \ml{\h,r}{X}$ (or equivalently, there exists a left $\h$-multiplications on $X$ making it into a Hilbert $\h$-bimodule). 

\smnoind
(c) Suppose that $\Lambda_1$ and $\Lambda_2$ are two unital $\R$-$*$-homomorphisms from $\h$ to $\ml{\h,r}{X}$.
Then there exists a unitary $U\in \ml{\h,r}{X}$ such that $\Lambda_2 = Ad(U)\circ \Lambda_1$. 
Consequently, if $Y_1$ and $Y_2$ are two Hilbert $\h$-bimodules such that they are isomorphic as Hilbert right $\h$-modules, then they are isomorphic as Hilbert $\h$-bimodules.
\end{prop}
\begin{prf}
(a) Suppose that $E$ is a closed convex subset of $X$ and $x\in X$. 
By the argument as in the case of complex Hilbert spaces, there exists a unique element $e\in E$ such that $\|x-e\| = {\rm dist}(x, E)$. 
Again, using the same argument as that for complex Hilbert spaces, for any closed $\h$-submodule $Y\subseteq X$ and any $x\notin Y$, there exists $z\in X\setminus \{0\}$ such that $\langle y, z \rangle = 0$ for any $y\in Y$. 
Thus, one can use a Zorn's lemma type argument to show that there exists a subset $\{e_i\}_{i\in I} \subseteq X$ such that $\langle e_i, e_j \rangle  = \delta_{i,j}$ and the right $\h$-linear span of $\{e_i\}_{i\in I}$ is dense in $X$. 
We now define on $X$ a real inner product $\re \langle \cdot, \cdot \rangle$ and let $H$ be the resulting real Hilbert space. 
Note that the norm defined by $\langle\cdot, \cdot \rangle$ and $\re \langle \cdot, \cdot \rangle$ are the same. 
If $D_1$ and $D_2$ are two orthonormal bases of $X$, then it is clear that $\{ xe: x\in D_1; e\in B\}$ and $\{ ye: y\in D_2; e\in B\}$ are two orthonormal bases of $H$. 
Consequently, $4\cdot {\rm card}(D_1) = 4\cdot {\rm card}(D_2)$. 

\smnoind
(b) Let $\{ e_i \}_{i\in I}$ be an orthonormal basis of $X$. 
Suppose that 
$$K_0 \ = \ \{ (\lambda_i)_{i\in I} \in l^2(I; \R):  \lambda_i = 0 {\rm\ except\ for\ finite\ numbers\ of}\ i\}.$$ 
Then the map $\Omega: K_0\otimes \h \rightarrow X$ given by $\Omega((\lambda_i)_{i\in I} \otimes \alpha) = \sum_{i\in I} e_i \alpha$ is a well defined right $\h$-module map. 
Moreover, for any finite subset $F\subseteq I$ and any $\alpha_i, \beta_i\in \h$ ($i\in F$), we have 
$$\left\langle \Omega((\alpha_i)_{i\in I}),  \Omega((\beta_i)_{i\in I}) \right\rangle \ = \ \sum_{i\in F} \alpha_i^*\beta_i$$
(we set $\alpha_i = 0 = \beta_i$ if $i\notin F$). 
Therefore, $\Omega$ is a Hilbert right $\h$-module isomorphism from $l^2(I; \R)\otimes \h$ (which is naturally a Hilbert $\h$-bimodule) onto $X$. 

\smnoind
(c) Suppose that $Y_1$ and $Y_2$ are the resulting Hilbert $\h$-bimodules defined by $\Lambda_1$ and $\Lambda_2$ respectively. 
Then for any $x,y\in (Y_i)_\real$ ($i=1,2$) and $\alpha\in \h$, 
$$\alpha \langle x,y \rangle 
\ = \ \langle x\alpha^*,y \rangle 
\ = \ \langle \alpha^* x ,y \rangle 
\ = \ \langle  x , \alpha y \rangle 
\ = \ \langle  x , y \alpha \rangle 
\ = \ \langle  x , y \rangle \alpha$$
which implies that $\langle x,y \rangle \in \R$. 
Thus, $(Y_i)_\real$ is a $\R$-Hilbert space and it is not hard to check that $Y_i \cong (Y_i)_\real \otimes \h$ as Hilbert $\h$-bimodules ($i=1,2$).  
If $E_i$ is an orthonormal basis for $Y_i$ ($i=1,2$), then $E_1$ and $E_2$ are both orthonormal bases for $X$ and hence are of the same cardinality (by part (a)). 
Suppose that $U\in \ml{\h,r}{X}$ is given by the bijection between $E_1$ and $E_2$. 
Then clearly $U$ is a unitary and $\Lambda_2 = Ad(U)\circ \Lambda_1$. 
The last statement is clear. 
\end{prf}

\medskip

Since any Hilbert right $\h$-module is automatically a Hilbert $\h$-bimodule in an essentially unique way, we believe that one can use the method in this paper to give easier arguments of the main results in \cite{Tor} and \cite{Vis}. 

\medskip

We have noted in the above that there are many different ``norm-quaternionizations'' for a given $\R$-normed space.  
However, the ``quaternionization'' of a $\R$-Hilbert space is unique in the sense that there is only one ``norm-quaternionizations'' that comes from a ``quaternion Hilbert space''. 

\medskip

\begin{thm}
\label{hips=iphb}
The following four categories are equivalent:

\smnoind
1. the category $\mathcal{A}$ of two-sided $\h$-Hilbert spaces; 

\smnoind
2. the category $\mathcal{B}$ of Hilbert $\h$-bimodules;

\smnoind
3. the category $\mathcal{C}$ of Hilbert right $\h$-modules;

\smnoind
4. the category $\mathcal{D}$ of $\R$-Hilbert spaces.
\end{thm}
\begin{prf}
Firstly of all, note that one can use Proposition \ref{cp-hil} and its argument to show that there is a bijective correspondences amongst objects of $\mathcal{B}$, $\mathcal{C}$ and $\mathcal{D}$. 
Suppose that $(Y, \lala \cdot, \cdot \rara)$ is a two-sided $\h$-Hilbert space. 
Define $\langle \cdot , \cdot \rangle: Y\times Y \rightarrow \h$ by 
$$\langle x,y \rangle\ :=\ m (\lala x, y \rara) \qquad (x,y\in Y)$$
(where $m$ is the multiplication on $\h$).
Then it is clear that $\langle \cdot , \cdot \rangle$ is a $\h$-valued pre-inner product on $Y$. 
Moreover, if $\langle x , x \rangle = 0$ and $\lala x, x \rara = \sum_{i=1}^n \beta_i^*\otimes \beta_i$, then 
$\sum_{i=1}^n \beta_i^* \beta_i = \langle x , x \rangle = 0$ which implies that $\beta_i = 0$ (for $i=1,...,n$) and so $\lala x, x \rara =0$ or equivalently, $x = 0$.  
On the other hand, for any $x,y\in Y$, if $\lala x, y \rara = \sum_{i=1}^n \alpha_i \otimes \beta_i$, then 
for any $\gamma\in \h$, 
\begin{equation*}
\langle \gamma x, y \rangle \ = \ m\big(\lala x, y \rara (\gamma^*\otimes 1)\big) \ = \ m\left(\sum_{i=1}^n \alpha_i\gamma^*\otimes \beta_i\right) \ = \ m\left(\sum_{i=1}^n \alpha_i\otimes \gamma^*\beta_i\right) \ = \ \langle x, \gamma^* y \rangle.
\end{equation*}
Therefore, $Y$ becomes a Hilbert $\h$-bimodule under $\langle \cdot, \cdot \rangle$. 
Conversely, let $(Y, \langle \cdot, \cdot \rangle)$ be a Hilbert $\h$-bimodule. 
Then the argument of Proposition \ref{cp-hil}(c) shows that $Y_\real$ is a $\R$-Hilbert space under the restriction of the inner product $\langle \cdot, \cdot \rangle$ and one can easily check that the following defines a two-sided $\h$-inner product on $Y = Y_\real\otimes \h$ (see Example \ref{eg-hil}): 
$$\lala y\otimes \alpha, z \otimes \beta \rara \ := \ \langle y, z \rangle \ \alpha^*\otimes \beta \qquad (x,y\in Z; \alpha, \beta\in \h).$$ 
It is not hard to see that the above gives a bijective correspondence between the objects of $\cal A$ and those of $\cal B$. 
Finally, it is easily seen that the morphisms in all these categories are also in bijective correspondences. 
\end{prf}

\medskip

By Proposition \ref{cp-hil}, for any Hilbert $\h$-bimodule $K$, there exists a $\R$-Hilbert space $(H,(\cdot,\cdot))$ such that $K \cong H \otimes \h$.
Let $K^{op}$ be the $\R$-vector space $H \otimes \h$ equipped with new $\h$-multiplications and a new inner product:
$$ \alpha\bullet (x\otimes \gamma) \bullet \beta \ :=\ x\otimes \beta^* \gamma \alpha^* \quad {\rm and} \quad 
\la x\otimes \alpha, y\otimes \beta \ra \ := \ (y, x)\ \alpha\beta^* \qquad (\alpha,\beta,\gamma\in\h; x,y\in H).$$
Then $K^{op}$ is a Hilbert $\h$-bimodule and is called the \emph{opposite Hilbert $\h$-bimodule} of $K$.

\medskip

The following is another equivalent form of ``quaternion Hilbert spaces'' which is virtually weaker than Hilbert $\h$-bimodules. 
The idea of this result comes from the argument of \cite[Theorem 4]{Kul}. 

\medskip

\begin{prop}
\label{sh>q}
Let $(H, (\cdot, \cdot))$ be a $\R$-Hilbert space.  
The existence of a $*$-homomorphism $\pi: \h \rightarrow \cl_{\R}({H})$ is equivalent to the existence of a $\R$-Hilbert space $H_0$ such that $H \cong H_0 \otimes_{\rm Hil} \h$. 
In this case, there is a Hilbert $\h$-bimodule structure on $H$ (with $\h$-valued inner product $\la \cdot, \cdot\ra$) such that $(\cdot, \cdot) = \re\la \cdot, \cdot\ra$ and $\pi$ is induced by the left $\h$-multiplication. 
\end{prop}
\begin{prf}
Suppose that such a $H_0$ exists. 
If $\Lambda: \h \rightarrow \cl_\R(\h)$ is the $*$-homomorphism given by the left $\h$-multiplication, then $\pi := 1\otimes \Lambda: \h \rightarrow \cl_\R(H_0 \otimes_{\rm Hil} \h)$ is a $*$-homomorphism. 
Conversely, suppose such a $\pi$ exists and we define a right $\h$-multiplication and a $\h$-valued inner product on $H$ by 
\begin{equation}
\label{hrhm}
x\cdot \alpha := \pi(\alpha^*)(x) \qquad {\rm and} \qquad [x, y] := \sum_{e\in B} (x, \pi(e)(y))\ e \qquad (x,y\in H; \alpha\in \h).
\end{equation}
It is clear that $[ \cdot, \cdot ]$ is $\R$-bilinear. 
For any $x,y\in H$ and $f\in B$, 
$$[x, y \cdot f] 
\ = \ \sum_{e\in B} (x, \pi(ef^*)(y))\ e
\ = \ \sum_{g\in B} (x, \pi(g)(y))\ gf
\ = \ [x, y]\cdot f$$
and 
$$[y, x] 
\ = \ \sum_{e\in B} (\pi(e^*)(y), x)\ e
\ = \ \sum_{g\in B} (x, \pi(g)(y))\ g^*
\ = \ [x, y]^*.$$
Moreover, for any $x\in H$ and $e\in B\setminus \{1\}$, 
$$(x, \pi(e)(x))
\ = \ (\pi(e^*)(x), x)
\ = \ - (\pi(e)(x), x) 
\ = \ - (x, \pi(e)(x))$$ 
and so $[x, x] = (x, x)$. 
Consequently, $H$ is a Hilbert right $\h$-module, denoted by $K$, under the structure in (\ref{hrhm}). 
By Proposition \ref{cp-hil}, there is a $\R$-Hilbert space $H_0$ such that $K\cong H_0\otimes \h$ and so $H \cong H_0\otimes_{\rm Hil} \h$. 
Finally, there exists a left $\h$-multiplication on $K$ turning it into a Hilbert $\h$-bimodule (by Proposition \ref{cp-hil} again), also denoted by $K$. 
Now it is not hard to check that $\pi$ is given by the left multiplication on $K^{op}$ and $(\cdot, \cdot) = \re\la \cdot, \cdot\ra$ (where $\la \cdot, \cdot\ra$ is the $\h$-valued inner product of $K^{op}$). 
\end{prf}

\medskip

The equivalences of categories in the above tells us that for any inner product $\h$-bimodules $X$ and $Y$, the two vector spaces $\mcl{\h}{X}{Y}$ and $\mcl{\R}{X_\real}{Y_\real}$ are isomorphic. 
In fact, they are isometrically isomorphic. 

\medskip

\begin{prop}
\label{mor-inn-pd-sp}
Let $X$ and $Y$ be inner product $\h$-bimodules. 

\smnoind
(a) $X \cong X_\real \otimes_{\rm Hil} \h$ as normed spaces (where $\otimes_{\rm Hil}$ is the Hilbert space tensor product) under the canonical $\h$-bimodule isomorphism $\Delta$ that sends $\sum_{e\in B} x_e e$ to $\sum_{e\in B} x_e \otimes e$ (where $x_e\in X_\real$). 

\smnoind
(b) The canonical isomorphism $\Phi: \mcl{\h}{X}{Y} \rightarrow \mcl{\R}{X_\real}{Y_\real}$ is isometric. 
\end{prop}
\begin{prf}
(a) Note that if $e,f\in B$ such that $e^*f \neq 1$, then $f^*e = - e^*f$. 
Thus, 
$$\left\|\sum_{e\in B} x_e e \right\|^2 \ =\ \left\|\sum_{e\in B} \langle x_e, x_e\rangle \right\| \ =\ \sum_{e\in B} \|x_e\|^2 \ = \ \|\Delta(x)\|^2.$$

\smnoind
(b) As $\Phi(T) = T\!\mid\!_{X_\real}$, it is clear that $\|\Phi(T)\| \leq \|T\|$. 
For any $x= \sum_{e\in B} x_ee\in X$ ($x_e\in X_\real$), we have, using the same argument as in part (a), 
$$\|T(x)\|^2 
\ = \ \left\|\sum_{e\in B} \Phi(T)(x_e)e\right\|^2 
\ = \ \sum_{e\in B} \|\Phi(T)(x_e)\|^2
\ \leq \ \|\Phi(T)\|^2 \sum_{e\in B} \|x_e\|^2
\ = \ \|\Phi(T)\|^2\|x\|^2. $$
\end{prf}

\medskip

\begin{rem}
\label{rem-quat}
\rm
(a) Note that the category of two-sided $\h$-inner product spaces, the category of inner product $\h$-bimodules and the category of $\R$-inner product spaces are also equivalent but we do not know if they are equivalent to inner product right $\h$-modules.  

\smnoind
(b) For any two $\h$-normed spaces $X$ and $Y$, the two Banach spaces
$\mcl{\h}{X}{Y}$ and $\mcl{\R}{X_\real}{Y_\real}$ are isomorphic (because $\| T\!\!\mid\!\!_{X_\real} \| \leq \|T\| \leq 4 \ \| T\!\!\mid\!\!_{X_\real} \|$ for any $T\in \mcl{\h}{X}{Y}$) but they are in general not isometrically isomorphic.

\smnoind
(c) By the argument of Theorem \ref{hips=iphb}, if $(Y, \lala \cdot, \cdot \rara)$ is a two-sided $\h$-Hilbert space and $x,y\in Y$, then $\|m(\lala x, y \rara)\| \leq \|x\|\ \|y\|$ (because of the corresponding result for Hilbert right $\h$-module which can be shown using a similar argument as that for Hilbert modules over complex $C^*$-algebras as in \cite{Lan}). 
\end{rem}

\medskip

Our next task is the Riesz representation type theorem for ``quaternion Hilbert spaces''. 
Let us first consider the notion of the dual space of a $\h$-normed space. 
If $X$ is a $\h$-normed space, then the space of quaternion linear functional $\mcl{\h}{X}{\h}$ is in general not a $\h$-vector space.
For example, when $X = \h$, we have $\mcl{\h}{X}{\h} = \R$. 

\medskip

We recall from \cite{Ng-dual-quat} the following definition of the dual object of $X$ (which is motivated by the duality of operator spaces; see e.g. \cite{Ng-reg}):
$$X^r \ := \ \mcl{\h}{X}{\hthr{\ep}}$$
where $\hthr{\ep}$ is the $\h$-bimodule $\hthr{\!}$ equipped with the injective tensor norm (see e.g. \cite{Ryan} for a definition).
$X^r$ is a $\h$-Banach space under the multiplication: 
$$(\alpha \cdot T \cdot \beta)(x)\ =\ (\alpha \otimes 1) T(x) (\beta\otimes 1)\qquad \qquad (\alpha,\beta\in \h; x\in X; T\in X^r).$$ 
Our next remark shows that $X^r$ is a nice duality because 
$$\mcl{\h}{X}{\hthr{\ep}} \ \cong \ \mcl{\R}{X_\real}{\h}$$
as $\h$-Banach spaces (compare with Lemma \ref{lh=lr}(b)). 

\medskip

\begin{rem}
\label{h-hh=r-h}
\rm
(a) As $\h$ is finite dimensional, $\mcl{\h}{X}{\hthr{\!}} \cong \h \otimes \mcl{\h}{X}{\h}$ as vector spaces. 
Moreover, for any $T\in \mcl{\h}{X}{\hthr{\ep}}$, 
$$\|T\|_{\mathcal{L}_{\h}(X;\hthr{\ep})} 
\ = \ \sup_{\|x\|\leq 1} \sup_{\|f\|\leq 1} \|(f\otimes \id)T(x)\|
\ = \ \sup_{\|f\|\leq 1} \|(f\otimes \id)\circ T\|_{\mathcal{L}_{\h}(X;\h)}
\ = \ \|T\|_{\h \otimes_\ep \mathcal{L}_{\h}(X;\h)}.$$
Therefore, by Lemma \ref{lh=lr}(b), we have  
$$X^r 
\ = \ \mcl{\h}{X}{\hthr{\ep}} 
\ \cong \ \h \otimes_\epsilon \mcl{\h}{X}{\h}
\ \cong \ \h \otimes_\epsilon X_\real^* 
\ \cong \ \mcl{\R}{X_\real}{\h}$$ 
as $\h$-Banach spaces (the last equivalence follows from a similar argument as the first one)
and the isometry from $\mcl{\h}{X}{\hthr{\ep}}$ to $\mcl{\R}{X_\real}{\h}$ is the map $\Phi$ in Proposition \ref{mor-inn-pd-sp}(b). 

\smnoind
(b) Let $Y$ be a inner product $\h$-bimodule (or equivalently, a two-sided $\h$-inner product space).
By part (a) and Proposition \ref{mor-inn-pd-sp}(b), the identity map 
$${\rm Id}:\mcl{\h}{Y}{\hthr{\ep}}\rightarrow \mcl{\h}{Y}{\hthr{\rm Hil}}$$ 
is an isometry. 
Note, however, that the norms on $\h\otimes_\ep \h$ and $\h \otimes_{\rm Hil} \h$ are not the same, e.g. 
if $\theta = \sum_{e\in B} e\otimes e$, then $\|\theta\|_{\rm Hil} = 2$ but $\|\theta\|_\ep = 1$ 
(because under the canonical isomorphism $\h \otimes_\ep \h \cong \h^* \otimes_\ep \h \cong \ml{\R}{\h}$, the corresponding element of $\theta$ is the identity map in $\ml{\R}{\h}$; see the identification of $\h^*\cong \h$ in the paragraph preceding Definition \ref{def-q-norm}). 
\end{rem}

\medskip

\begin{lem}
\label{rem-inn-pd}
(a) If $X$ is a $\h$-normed space and $T\in X^r$, we set 
$$\|T\|_L \ :=\ \|m\circ T\|$$
where $m$ is the multiplication on $\h$. 
Then $\|T \| \leq \|T\|_L$ ($T\in X^r$) and $\|\cdot\|_L$ is a left $\h$-module norm on $X^r$. 

\smnoind
(b) Let $(Y, \lala \cdot, \cdot \rara)$ be a two-sided $\h$-inner product space. 
For any $y\in Y$, we can define $T_y \in L_{\h}(Y;\hthr{\!})$ by 
$$T_y(x)\ :=\ \lala y, x \rara \qquad (x\in Y).$$
Then $T_y \in Y^r$ and $\|T_y\| \leq \|y\|$. 
\end{lem}
\begin{prf}
(a) By Remark \ref{h-hh=r-h}(a), for any $T\in X^r$, 
$$\|T\|\ = \ \|\Phi(T)\| \ = \ \big\|(m\circ T)\!\mid\!_{X_\real}\big\| \ \leq \ \|T\|_L.$$
Consequently, $\|\cdot\|_L$ is a norm and it is clear that $\|\cdot\|_L$ is a left $\h$-module norm according to the scalar multiplication defined on $X^r$. 

\smnoind
(b) It is clear that $T_y \in L_{\h}(Y;\hthr{\!})$ and $\|T_y(x)\| = \|m(\lala y, x \rara)\| \leq \|y\|\ \|x\|$ for any $x\in Y$ (see Remark \ref{rem-quat}(c)). 
\end{prf}

\medskip

The following example shows that in general we do not have $\|T\| = \|T\|_L$. 

\medskip

\begin{eg}
Let $X$ be the Hilbert $\h$-bimodule $\mathbb{C} \otimes \h$ and $T\in X^r$ be defined by 
$T(1 \otimes \beta_1 + i \otimes \beta_i) := J(1)\otimes \beta_1  + J(i) \otimes \beta_i$ for any $\beta_1, \beta_i \in \h$ (where $J$ is the ``forgettable isometry'' from $\mathbb{C}$ to $\h$). 
Then Remark \ref{h-hh=r-h}(a) tells us that $\| T \| = \|J\| = 1$. 
On the other hand, 
$$\| T \|_L \ = \ \sup\{ \|\beta_1 + i\beta_i \|: \beta_1, \beta_i \in \h; \|1 \otimes \beta_1 + i \otimes \beta_i\| \leq 1 \}.$$
By considering $\beta_1 = i/\sqrt{2}$ and $\beta_i = 1/\sqrt{2}$, we see that $\| T \|_L \geq \sqrt{2}$. 
\end{eg}

\medskip

It follows directly from Theorem \ref{hips=iphb} that any element in $\mcl{\h}{Y}{\h}$ comes from an element in $Y_\real$. 
The following result gives a complete version for a Riesz representation type theorem (which reproduces all the elements in $Y$). 

\medskip

\begin{thm}(Riesz Representation theorem)
\label{rrt}
Suppose that $Y$ is a Hilbert $\h$-bimodule. 
Then for any $T\in Y^r$, there exists a unique $y\in Y$ such that 
$$m(T(x))\ = \ \langle y, x\rangle \qquad \qquad (x\in X)$$
and $\|T\|_L\ = \ \|y\|$. 
\end{thm}
\begin{prf}
Since $Y^r \cong Y_\real^*\otimes_\ep \h \cong Y_\real \otimes_\ep \h$ as Banach spaces (note that $Y_\real$ is a $\R$-Hilbert space because of Theorem \ref{hips=iphb}), there exist $y_e\in Y_\real$ ($e\in B$) such that for any $z_e\in Y_\real$, we have 
$$T\left(\sum_{e\in B} z_e e\right)\ =\ \sum_{e,f\in B} \lala y_f,z_e \rara (f\otimes e)$$ 
where $\lala \cdot, \cdot \rara$ is the corresponding two-sided $\h$-inner product given by Theorem \ref{hips=iphb} and so $y = \sum_{e\in B} e^* y_e$ will satisfy the first equality in the statement. 
The second equality follows directly from the definition of $\|\cdot \|_L$. 

\end{prf}

\bigskip

\bigskip

\section{Quaternion $B^*$-algebras}

\bigskip

In \cite{HS} and \cite{Kul}, (unital) real $B^*$-algebras that contain (unital) subalgebras $*$-isomorphic to $\h$ are considered. 
This kind of algebras are the same as quaternion $B^*$-algebras as defined in the following (see the discussion in Remark \ref{unital} and Remark \ref{rem-h-bs}(c) below). 
Using our observation concerning the real parts of quaternion vector spaces, we can obtain easily the improved versions of the main results in \cite{Kul}. 
Let us start with the definition of quaternion algebras. 

\medskip

\begin{defn}
\label{def-bsalg}
(a) Let $A$ be a $\h$-vector space as well as a $\R$-algebra (with respect to the induced $\R$-vector space structure on $A$). 
Then $A$ is called a \emph{$\h$-algebra} if the following conditions are satisfied for any $a,b\in A$ and $\gamma\in \h$:
$$(\gamma a)b\ = \ \gamma (ab), \qquad (a\gamma)b\ =\ a(\gamma b) \qquad {\rm and} \qquad (ab)\gamma = a (b\gamma).$$
Moreover, a $\R$-involution $^*$ on $A$ is called a \emph{$\h$-involution} if for any $\alpha \in \h$ and $b\in A$, 
$$(\alpha b)^*\ =\ b^*\alpha^* \qquad {\rm and} \qquad (b\alpha)^* \ = \ \alpha^*b^*.$$ 
In this case, $A$ is called a \emph{$\h$-involutive algebra}. 

\smnoind
(b) Suppose that $A$ is a $\h$-involutive algebra with a $\h$-norm $\|\cdot\|$. 
Then $A$ is called a \emph{$\h$-$B^*$-algebra} if it is complete under $\|\cdot\|$ and $\|\cdot\|$ satisfies the following properties:  
$$\|ab\|\ \leq \ \|a\| \ \|b\| \qquad {\rm and} \qquad \|a^*a\| = \|a\|^2 \qquad \qquad (a,b\in A).$$
\end{defn}

\medskip

\begin{rem}
\label{rem-h-bs}
\rm
(a) We recall that a $\R$-Banach-$*$-algebra is a \emph{$\R$-$B^*$-algebra} if its norm satisfies the second equality in Definition \ref{def-bsalg}(b). 
Moreover, a $\R$-$B^*$-algebra $A$ is said to be \emph{hermitian} if $\sigma_{A+iA}(a) \subseteq \R$ for any element $a\in A$ with $a^* = a$. 
Note that whether $A$ is hermitian depends only on the algebraic structure of $A$ (i.e. unrelated to the norm on $A$). 

\smnoind
(b) A \emph{$\R$-$C^*$-algebra} is a closed $*$-subalgebra of $\ml{\R}{H}$ for a $\R$-Hilbert space $H$. 
For a $\R$-$B^*$-algebra $B$, the following statements are equivalent. 
\begin{enumerate}
\item[(i).] $B$ is a $\R$-$C^*$-algebra. 

\item[(ii).] There exists an injective $*$-representation $\pi$ of $B$ on a $\R$-Hilbert space such that $\pi(B)$ is closed. 

\item[(iii).] $B$ is hermitian. 

\item[(vi).] $B$ can be ``complexified'' into a $\C$-$C^*$-algebra
\end{enumerate}
In fact, the equivalences of statements (i), (iii) and (iv) can be found in \cite[5.1.2 \& 5.2.11]{Li}. 
On the other hand, it is clear that statement (i) implies statement (ii). 
Moreover, if statement (ii) holds, then $\pi(B)$ is hermitian (as it is a $\R$-$C^*$-algebra) and so is $B$ (because $\pi: B \rightarrow \pi(B)$ is a bijective $*$-homomorphism). 

\smnoind
(c) If $A$ is a unital $\R$-$B^*$-algebra with an injective unital $*$-homomorphism $\phi: \h \rightarrow A$, then $A$ is a $\h$-$B^*$-algebra under the $\h$-scalar multiplication induced by $\phi$. 

\smnoind
(d) If $A$ and $B$ are $\R$-$C^*$-algebras, then any injective $*$-homomorphism $\varphi: A \rightarrow B$ is isometric. 
In fact, as in part (b), for any $\R$-$C^*$-algebra $A$, there exists a $\C$-$C^*$-algebra $A_C$ such that $A$ is a closed $\R$-$*$-subalgebra of $A_C$ and $A_C = A + Ai$. 
Since $A_C \cong A \otimes \C$ as complex $*$-algebras, $\varphi$ can be extended to an injective (and hence isometric) complex $*$-homomorphism from $A_C$ to $B_C$. 
\end{rem}

\medskip

\begin{eg}
\label{eg-h-bs}
Suppose that $(K, \langle \cdot, \cdot \rangle)$ is a Hilbert $\h$-bimodule. 

\smnoind
(a) By Remark \ref{rem-hil}(d), the left scalar multiplication on $K$ defines a $\R$-$*$-homomorphism $\Theta_K: \h \rightarrow \ml{\h,r}{K}$. 
Therefore, $\ml{\h,r}{K}$ is canonically a $\h$-$B^*$-algebra. 
Note that there is a canonical injective unital $*$-homomorphism 
$$J_K: \ml{\R}{K_\real}\otimes \h \rightarrow \ml{\h,r}{K}$$ 
given by $J_K(T\otimes \alpha)(x\otimes \beta) = T(x) \otimes \alpha\beta$ ($T\in \ml{\R}{K_\real}$, $\alpha,\beta \in \h$ and $x\in K_\real$). 
In this case, we have $\Theta_K = J_K \circ (1\otimes \Theta_\h)$.  

\smnoind
(b) Since $\ml{\h,r}{K}\subseteq \ml{\R}{K}$ as unital $\R$-$B^*$-algebra, we see that $\ml{\R}{K}$ is also a $\h$-$B^*$-algebra canonically. 
\end{eg}

\medskip

The proof of the following lemma is clear and will be omitted. 

\medskip

\begin{lem}
\label{dec-h-alg}
If $A$ is a $\h$-$B^*$-algebra, then $A_\real$ is a $\R$-$B^*$-subalgebra of $A$ and the map 
$\varphi: A\rightarrow A_\real\otimes \h$ given by 
$$\varphi: a\ \longmapsto\ \sum_{e\in B} \re( e^* a) \otimes e \qquad (a\in A)$$ 
is a $\h$-$*$-isomorphism.
\end{lem}

\medskip

\begin{thm}(Gelfand-Naimark theorem)
\label{gn-thm}
For any $\h$-$B^*$-algebra $A$, there exists a Hilbert $\h$-bimodule $H$ and an isometric $\h$-$*$-homomorphism from $A$ to $\ml{\h,r}{H}$. 
\end{thm}
\begin{prf}
Lemma \ref{dec-h-alg} tells us that $A_\real$ is a $\R$-$B^*$-algebra. 
It is easy to check that $A_{1,i} = A_\real + i A_\real$ (see Proposition \ref{complex}) is a $\mathbb{C}$-$B^*$-algebra (equivalently, a $\mathbb{C}$-$C^*$-algebra). 
An isometric $*$-representation of $A_{1,i}$ on a $\mathbb{C}$-Hilbert space $(K, (\cdot, \cdot))$ (which always exists) induces an isometric $*$-representation $\pi$ of $A_\real$ on the $\R$-Hilbert space $K_0 = (K, \re(\cdot, \cdot))$. 
Therefore, $\pi \otimes \Theta_\h$ (where $\Theta_\h$ is the map as in Example \ref{eg-h-bs}(a)) is an injective $*$-representation of $A\cong A_\real \otimes \h$ on the $\R$-Hilbert space $K_0\otimes l^{(4)}_2$ and it is not hard to see that $(\pi \otimes \Theta_\h)(A)$ is closed in $\ml{\R}{K_0\otimes \h}$. 
Consequently, $A$ is actually a $\R$-$C^*$-algebra (see Remark \ref{rem-h-bs}(b)). 
On the other hand, if $H = K_0 \otimes \h$, then $J_H \circ (\pi \otimes \id_\h)$ (where $J_H$ is the map in Example \ref{eg-h-bs}(a)) induced an injective $\h$-$*$-homomorphism from $A\cong A_\real \otimes \h$ to $\ml{\h,r}{H}$. 
This $*$-homomorphism is isometric because of Remark \ref{rem-h-bs}(d). 
\end{prf}

\medskip

Using this theorem as well as the last statement of Example \ref{eg-h-bs}(b), we can obtain the following extension of the main result in \cite{Kul}. 

\medskip

\begin{cor}(Alternative form of the Gelfand-Naimark theorem)
\label{imp-kul2}
For any $\h$-$B^*$-algebra $A$, there exists a Hilbert $\h$-bimodule $H$ and an isometric $\h$-$*$-homomorphism from $A$ to $\ml{\R}{H}$. 
\end{cor}

\medskip

Note that a unital $\R$-$B^*$-algebra containing a unital subalgebra $*$-isomorphic to $\h$ will automatically be a $\h$-$B^*$-algebra and so is a subalgebra of some $\ml{\R}{H}$ by the above corollary. 

\medskip

\begin{cor}
Suppose that $A$ is a unital $\R$-$B^*$-algebra. 
If there exists a unital $*$-homomophism from $\h$ to $A$, then $A$ is hermitian. 
\end{cor}

\medskip

\begin{rem}
\label{rem-kul2}
\rm
(a) As stated in Remark \ref{rem-hil}, because the quaternionic Hilbert space as defined in \cite[Definition 2]{Kul} is automatically unital, one needs to assume that the real $B^*$-algebra in \cite[Theorem 4]{Kul} is unital and contains a unital $*$-subalgebra $*$-isomorphic to $\h$. 
Therefore, \cite[Theorem 4]{Kul} is a particular case of Corollary \ref{imp-kul2}. 
Note that we do not need to assume, a prior, that $A$ is hermitian (as this condition automatically holds for $\h$-$B^*$-algebras). 

\smnoind
(b) Another interpretation of Theorem \ref{gn-thm} is that $\h$-$B^*$-algebras and \emph{$\h$-$C^*$-algebras} (which is defined as closed $\h$-$*$-subalgebras of $\ml{\h,r}{H}$ for some Hilbert $\h$-bimodule $H$) are the same.

\smnoind
(c) An application of Theorem \ref{gn-thm} is the following. 
If one want to define quaternionic operator spaces, the first problem is whether one should consider subspaces of $\ml{\R}{H}$ (where $H$ is a Hilbert $\h$-bimodule) or $\ml{\h,r}{K}$ (where $K$ is a Hilbert right $\h$-module). 
Theorems \ref{hips=iphb} and \ref{gn-thm} tell us that these two notions are actually the same. 
\end{rem}

\medskip

It is clear that if $A$ is a $\h$-$B^*$-algebra, then $A_\real$ is a $\R$-$B^*$-algebra and Theorem \ref{gn-thm} shows that $A_\real$ is a $\R$-$C^*$-algebra.
Conversely, if $B$ is a $\R$-$C^*$-algebra, then there exists a norm on $B\otimes \h$ which turns it into a $\h$-$B^*$-algebra (see again the argument of Theorem \ref{gn-thm}). 
Remark \ref{rem-h-bs}(d) tells us that such a norm on $B\otimes \h$ is unique. 
Furthermore, it is also clear that $\R$-$*$-homomorphisms between $\R$-$C^*$-algebras correspond bijectively to $\h$-$*$-homomorphisms between their tensor products with $\h$. 
Thus, we have the following corollary. 

\medskip

\begin{cor}
The category of $\h$-$B^*$-algebras is equivalent to the category of $\R$-$C^*$-algebras.
\end{cor}

\medskip

We can also obtained very easily the following improvement of \cite[Corollary 6]{Kul} (note that the proof in \cite{Kul} required the main result of \cite{Kul0} and the statement in \cite[Corollary 6]{Kul} includes the hermitian assumption). 

\medskip

\begin{thm}(Gelfand theorem)
\label{imp-kul1}
Let $A$ be a $\h$-$B^*$-algebra. 
Then $A \cong C_0(\Omega;\h)$ (as $\h$-$B^*$-algebra) for a locally compact Hausdorff space $\Omega$ if and only if all elements in $A$ are normal. 
\end{thm}
\begin{prf}
It is easy to check that if $A \cong C_0(\Omega;\h)$, then any element in $A$ is normal. 
Conversely, suppose that all elements in $A$ are normal. 
Take any $a,b,c,d\in A_\real$. 
From the equalities $(a^*+b^*)(a+b) = (a+b)(a^*+b^*)$ and $(a^*-b^*i)(a+bi) = (a+bi)(a^*-b^*i)$,
we know that 
$$a^*b + b^*a = ba^* + ab^* \qquad {\rm and} \qquad a^*b - b^*a\ =\ ba^* - ab^*$$
and hence $a^*b = b a^*$ which implies that $A_\real$ is commutative. 
On the other hand, by comparing the ``$i$-coefficients'' of the two expressions: 
$(a+bi+cj+dk)^*(a+bi+cj+dk)$ and $(a+bi+cj+dk)(a+bi+cj+dk)^*$, we see that 
$$a^*b - b^*a - c^*d + d^*c \ =\ -ab^* + ba^* - cd^* + dc^*$$ 
and so, $d^*c = dc^*$. 
Now, if we replace $c$ by an approximate identity for the $\R$-$C^*$-algebra $A_\real$ (see e.g. \cite[5.2.4]{Li}), we see that $d^* = d$ (i.e. the involution is trivial on $A_\real$). 
Since $A_{1,i} = A_\real + A_\real i$ is a commutative $\C$-$C^*$-algebra, there exists a locally compact Hausdorff space $\Omega$ such that $A_{1,i} \cong C_0(\Omega, \C)$. 
Because of the relation $A_\real = (A_{1,i})_{sa}$ (since the involution is trivial on $A_\real$), we see that $A_\real \cong C_0(\Omega, \R)$ and thus, $A \cong C_0(\Omega, \h)$. 
\end{prf}

\medskip

A $\h$-$B^*$-algebra $A$ is said to be \emph{normal} if every element in $A$ is normal. 

\medskip

\begin{cor}
The following categories are equivalent:

\smnoind
1. the category of compact Hausdorff spaces with continuous maps as morphisms; 

\smnoind
2. the category of unital normal $\h$-$B^*$-algebras with unital $\h$-$*$-homomorphisms as morphisms;

\smnoind
3. the category of unital commutative $\mathbb{C}$-$C^*$-algebras with unital $\mathbb{C}$-$*$-homomorphisms as morphisms.
\end{cor}

\bigskip

\noindent
School of Mathematical Sciences and LPMC,
Nankai University,
Tianjin 300071,
China.

\smallskip\noindent
E-mail address: ckng@nankai.edu.cn


\bigskip

\bigskip

\begin{thebibliography}{99}

\bibitem{AJ} M. Abel and K. Jarosz, Noncommutative uniform algebras, Studia Math. 162 (2004), 213-218. 

\bibitem{AK} S. Agrawal and S. H. Kulkarni, An analogue of the Riesz-representation theorem, Novi Sad J. Math. 30 (2000), 143-154.

\bibitem{AK2} S. Agrawal and S. H. Kulkarni, Dual spaces of quaternion normed linear spaces and reflexivity, J. Anal. 8 (2000), 79-90.

\bibitem{BCT} L. C. Biedenharn, G. Cassinelli and P. Truini, Imprimitivity theorem and quaternionic quantum mechanics, \emph{Proceedings of the Third Workshop on Lie-Admissible Formulations (Univ. Massachusetts, Boston, Mass., 1980), Part B}, Hadronic J. 4 (1980/81), 981-994.

\bibitem{CT} G. Cassinelli and P. Truini, Quantum mechanics of the quaternionic Hilbert spaces based upon the imprimitivity theorem, Rep. Math. Phys. 21 (1985), 43-64. 

\bibitem{BH} L. C. Biedenharn and L. P. Horwitz, Quaternion quantum mechanics: second quantization and gauge fields, Ann. Phys. 157 (1984), 432-488. 

\bibitem{HR-ten} L. P. Horwitz and A. Razon, Tensor product of quaternion Hilbert modules, Acta Appl. Math. 24 (1991), 141-178. 

\bibitem{HR-proj} L. P. Horwitz and A. Razon, Projection operators and states in the tensor product of quaternion Hilbert modules, Acta Appl. Math. 24 (1991), 179-194. 

\bibitem{HR-uniq} L. P. Horwitz and A. Razon, Uniqueness of the scalar product in the tensor product of quaternion Hilbert modules, J. Math. Phys. 33 (1992), 3098-3104. 

\bibitem{HS} L. P. Horwitz and A. Soffer, $B^*$-algebra representations in quaternionic Hilbert module, J. Math. Phys. 24 (1984), 2780-2782. 

\bibitem{Kul0} S. H. Kulkarni, Representation of a class of real $B^*$-algebras as algebras of quaternion-valued functions, Proc. Amer. Math. Soc. 116 (1992), 61-66. 

\bibitem{Kul} S. H. Kulkarni, Representation of a real $B^*$-algebra on a quaternionic Hilbert space, Proc. Amer. Math. Soc. 121 (1994), 505-509. 

\bibitem{Lan} E. C. Lance, ``Hilbert $C\sp *$-modules: A toolkit for operator algebraists'', Lond. Math. Soc. Lect. Note Ser. 210, Cambridge University Press, Cambridge (1995).

\bibitem{Li} B. R. Li, \emph{Real Operator algebras}, World Scientific Publishing Co., River Edge NJ (2003).

\bibitem{Lud} S. V. Ludkovsky, Unbounded operators on Banach spaces over the quaternion field, arXiv: math.OA/0404444 v1 (2004).

\bibitem{NV} S. Natarajan and K. Viswanath, Quaternionic representations of compact metric groups, 
J. Math. Phys. 8 (1967) 582-589.

\bibitem{Ng-reg} C. K. Ng, Regular normed bimodules, preprint. 

\bibitem{Ng-dual-quat} C. K. Ng, Quaternion Banach spaces, in preparation. 

\bibitem{Ryan} R. A. Ryan, \emph{Introduction to tensor products of Banach spaces}, Springer Monographs in Mathematics, Springer-Verlag, London (2002).

\bibitem{Suh} G. Suhomlinov, On extensions of linear functions in complex and quaternionic linear spaces, Mat. Sbornik 3 (1938), 353-358.

\bibitem{Tor0} A. Torga\v{s}ev, On the symmetric quaternionic Banach algebras I: Geljfand theory, Publ. Inst. Math. (Beograd) (N.S.) 24 (38) (1978), 173-188.

\bibitem{Tor} A. Torga\v{s}ev, Quaternionic operators with finite matrix trace, Int. Eq. Oper. Th. 23 (1995), 114-122. 

\bibitem{Tor2} A. Torga\v{s}ev, On reflexivity of a quaternion normed space, Novi Sad J. Math. 27 (1997), 57-64.

\bibitem{Tor-dual} A. Torga\v{s}ev, Dual space of a quaternion Hilbert space, Facta Univ. Ser. Math. Inform.  No. 14 (1999), 71-77. 

\bibitem{Vis} K. Viswanath, Normal operations on quaternionic Hilbert spaces, Trans. Amer. Math. Soc.  162 (1971), 337-350. 

\end{thebibliography}
\end{document}